\documentstyle[12pt]{amsart}

\newtheorem{thm}{Theorem}
\newtheorem{prop}{Proposition}

\newtheorem{cor}{Corollary}
\newtheorem{conj}{Conjecture}

\theoremstyle{remark}
\newtheorem{rem}{Remark}

\theoremstyle{definition}
\newtheorem{defn}{Definition}

\newcommand{\bZ}{\Bbb Z}

\newcommand{\bC}{\Bbb C}

\newcommand{\blowup}{\overline{\bC P^2}}

\title[Irreducible four--manifolds]
{On irreducible four--manifolds}
\author{D.~Kotschick}
\address{Mathematisches Institut, Universit\"at Basel,
Rheinsprung 21, 4051 Basel, Switzerland}

\begin{document}

\maketitle

\bigskip


\bigskip

\section{Introduction}\label{s:intro}

For many years, four--manifold folklore suggested that all
simply connected smooth four--manifolds should be connected
sums of complex algebraic surfaces, with both their complex
and non--complex orientations allowed\footnote{The $4$--sphere
is the empty connected sum.}. The first counterexamples
were constructed in 1990 by Gompf and Mrowka~\cite{GM}, and
many others followed. Then, Gompf~\cite{gompf} showed that many,
and possibly all, these counterexamples arise from symplectic
four--manifolds. Having no indication to the contrary, many
people have put forward the following:

\begin{conj}\label{sc}
Every smooth, closed, oriented and
simply connected $4$--manifold is the connected
sum of symplectic manifolds, with both the symplectic and the
opposite orientations allowed.
\end{conj}

\noindent
This conjecture is ambitious; it would imply the smooth
Poincar\'e conjecture.

Note that such a connected sum can have summands with
definite intersection forms, for example copies of $\bC P^2$.
In fact, any definite summand has a diagonalizable
intersection form by Donaldson's theorem~\cite{don},
and is therefore homeomorphic to $n\bC P^2$ by Freedman's
classification~\cite{freed}.

When the manifolds under consideration are not simply connected,
the situation is more complicated. Then there are obvious
counterexamples to Conjecture~\ref{sc}, e.g. rational homology
spheres
which are not homotopy spheres. Thus,
one has to allow definite summands which are more
general than $n\bC P^2$ or $n\blowup$. The natural conjecture is:

\begin{conj}\label{minimal}
Every smooth, closed and oriented $4$--manifold is the connected
sum of symplectic manifolds, with both the symplectic and the
opposite orientations allowed, and of some manifolds with
definite intersection forms.
\end{conj}

\noindent
This has occurred to several people, especially in the light
of the examples constructed in~\cite{KMT},
and has been dubbed the current ``minimal
conjecture'' by Taubes. In this note we show that it is
false\footnote{We do not have any proposal for a new
``minimal conjecture'' in the non--simply connected case.}.
Conjecture~\ref{sc} remains open.

Recall the following definition:

\begin{defn}
A smooth closed $4$--manifold $X$ is {\sl irreducible} if
for every smooth connected sum decomposition $X\cong X_1\#
X_2$ one of the summands $X_i$ must be a homotopy sphere.
\end{defn}

We will show:

\begin{thm}\label{t:counter}
There exist oriented irreducible $4$--manifolds $X$
with indefinite intersection forms and
with $\pi_1 (X)=\bZ_2$ and $b_2^+ (X)\equiv b_2^- (X)\equiv 0\pmod
2$.
\end{thm}

This follows from Proposition~\ref{p:irred1} in the next
section.

\begin{cor}
There exist orientable irreducible $4$--manifolds $X$
with indefinite intersection forms,
which are not almost complex (and therefore not complex
and not symplectic) with respect to either
orientation, and for which the Donaldson and Seiberg--Witten
invariants are not defined (or must vanish by definition).
\end{cor}

\noindent
If one drops the requirement that $X$ have indefinite
intersection form, rational homology spheres give
obvious examples.

Recall that a symplectic $4$--manifold is called minimal
if it contains no symplectically embedded $2$--sphere of
selfintersection $-1$. Conjectures~\ref{sc} and~\ref{minimal}
are complementary to Gompf's conjecture~\cite{gompf} that
minimal symplectic $4$--manifolds are irreducible. In
section~\ref{s:sympl} we deduce from the recent work of
Taubes~\cite{swg} on the Gromov and Seiberg--Witten invariants
that Gompf's conjecture is true in many cases, including all
simply connected manifolds with $b_2^+ > 1$. We shall prove:

\begin{thm}\label{t:sympl}
Let $X$ be a minimal symplectic $4$--manifold with $b_2^+ (X)>1$.
If $X\cong X_1\# X_2$ is a smooth connected sum decomposition
of $X$, then one of the $X_i$ is an integral homology sphere
whose fundamental group has no non--trivial finite quotient.
\end{thm}

\noindent
This strengthening of Proposition 1 in~\cite{KMT} goes a long
way towards confirming a conjecture made there.

\begin{cor}\label{c:minirr}
Minimal symplectic $4$--manifolds with $b_2^+ >1$
and with residually finite fundamental groups are irreducible.
\end{cor}

\noindent
See section~\ref{s:sympl} for a result in and comments
on the case $b_2^+=1$.

For K\"ahler surfaces, Theorem~\ref{t:sympl} and
Corollary~\ref{c:minirr} could be deduced easily from
the reduction of the Seiberg--Witten equation to
the K\"ahler vortex equation and the study of effective
divisors on complex surfaces, due to Kronheimer--Mrowka
and Witten~\cite{witten}.

\section{Irreducibility of quotient manifolds}\label{s:irred}

Theorem~\ref{t:counter} will follow from the following
application of the covering trick introduced in~\cite{IMRN}:

\begin{prop}\label{p:irred1}
Let $X$ be a smooth, closed, simply connected and oriented
spin $4$--manifold. If $b_2^+ (X)>1$, assume that $X$ has
a non--trivial Donaldson or
Seiberg--Witten invariant. Suppose a non--trivial finite group $G$
acts freely by orientation--preserving diffeomorphisms
of $X$. Then the quotient $Y=X/G$ is an orientable
irreducible $4$--manifold.
\end{prop}
\begin{pf}
Let $Y\cong M\# N$ be a smooth connected sum decomposition.
As $\pi_1 (Y)\cong G$ is finite, it cannot be a non--trivial
free product and we may assume $\pi_1 (M)\cong G$ and
$\pi_1 (N)\cong \{ 1\}$.

Let $d$ be the order of $G$. The connected sum decomposition
of $Y$ induces a connected sum decomposition $X\cong\overline{M}
\# dN$, where $\overline{M}$ is the universal covering of
$M$. As either $b_2^+ (X)\leq 1$ or $X$ is assumed to have
a non--trivial Donaldson
or Seiberg--Witten invariant, it follows that the intersection
form of $N$ is negative definite. By Donaldson's theorem~\cite{don}
it is diagonalizable over $\bZ$, and therefore either
trivial or odd.

On the other hand, the intersection form of $N$ must
be even, because it is a direct summand of the intersection
form of $X$, which is even because $X$ is spin. We conclude
$b_2 (Y)=0$. As $N$ is simply connected, it is a homotopy
sphere and $Y$ is irreducible.
\end{pf}

To obtain examples as in the statement of Theorem~\ref{t:counter},
take for $X$ the Fermat surface of degree $d\equiv 2 \pmod 4$
in $\bC P^3$ with $d\geq 6$. This is the surface defined in
homogeneous coordinates by
\begin{equation}\label{fermat}
x^d+y^d+z^d+t^d=0\ .
\end{equation}
It is simply connected, and spin because its canonical class
is the restriction of $(d-4)H$ to $X$, and therefore
$2$--divisible. Being an algebraic surface with $b_2^+>1$,
$X$ has non--trivial Donaldson
and Seiberg--Witten invariants. Furthermore,
the equation~\eqref{fermat} is invariant under complex conjugation
on $\bC P^3$, and has no non--trivial real solutions (because
$d$ is even). Thus,
complex conjugation acts freely on $X$ and the Proposition
shows that $Y=X/\bZ_2$ is irreducible. We have $\pi_1 (Y)=\bZ_2$.
Using the multiplicativity of the Euler characteristic and
the signature under finite unramified coverings, one
can calculate $b_2^{\pm}(Y)=\frac{1}{2}(b_2^{\pm}(X)-1)$ which
are both positive (because $d\geq 6$) and even (because
$d\equiv 2\bmod 4$).

This completes the proof of Theorem~\ref{t:counter}.

\begin{rem}
Wang~\cite{wang} has shown that the quotients of simply
connected minimal algebraic surfaces of general type by free
anti--holomorphic involutions have trivial Seiberg--Witten
invariants, even when the invariants do not vanish
by definition as in the above examples.
\end{rem}

\medskip

The assumptions in Proposition~\ref{p:irred1} are such that
$X$ has to be irreducible, and then $Y$ turns out irreducible
as well. Here is another result in the same spirit, but which
does not require a spin condition and uses instead
Corollary~\ref{c:minirr}. We could apply this to the examples
discussed above, but Proposition~\ref{p:irred1} is much more
elementary.

\begin{prop}\label{p:irred2}
Let $X$ be a closed, simply connected minimal symplectic
$4$--manifold with $b_2^+(X)>1$. Suppose a non--trivial
finite group $G$ acts freely by orientation--preserving
diffeomorphisms of $X$. Then the quotient $Y=X/G$ is an
orientable irreducible $4$--manifold.
\end{prop}
\begin{pf}
If $Y\cong M\# N$ with $\pi_1 (M)\cong G$ and $\pi_1 (N)
\cong \{ 1\}$, then $X\cong \overline{M}\# dN$.
Thus Corollary~\ref{c:minirr} implies that $N$ is
a homotopy sphere. Hence $Y$ is irreducible.
\end{pf}

\section{Irreducibility of symplectic manifolds}\label{s:sympl}

In this section we give the proof of Theorem~\ref{t:sympl}.
This requires some familiarity with Seiberg--Witten invariants,
particularly the work of Taubes~\cite{taubes,taubes2,swg}.
See also~\cite{KMT,KM,witten}.

Let $X$ be a closed symplectic $4$--manifold with $b_2^+ (X)>1$.
If $X$ splits as a connected sum $X\cong M\# N$, then by
Proposition 1 of~\cite{KMT} we may assume that $N$ has a
negative definite intersection form and that its fundamental group
has no non--trivial finite quotient. In particular $H_1 (N,\bZ )
=0$. This implies that the homology and cohomology of $N$
are torsion--free.

Donaldson's theorem about (non--simply connected) definite
manifolds~\cite{orient} implies that the intersection form
of $N$ is diagonalizable over $\bZ$. If $N$ is not an
integral homology sphere, let $e_1,\ldots ,e_n
\in H^2 (N,\bZ )$ be a basis with respect to which the
cup product form is the standard diagonal form. This basis
is unique up to permutations and sign changes.

It is a theorem of Taubes~\cite{taubes} that the
Seiberg--Witten invariants of $X$ are non--trivial
for the canonical $Spin^c$--structures with auxiliary
line bundles $\pm K_X$. Note that we can write
$$
K_X=K_M+\sum_{i=1}^n a_i e_i \ ,
$$
where $K_M\in H^2 (M,\bZ )$ and the $a_i$ are odd integers
because $a_i^2 = -1$ and $K_X$ is characteristic. Considering
$-K_X$ and using a family of Riemannian metrics which
stretches the neck connecting $M$ and $N$, we conclude
that $M$ has a non--trivial Seiberg--Witten invariant
for a $Spin^c$--structure with auxiliary line bundle $-K_M$.

Now we can reverse the process and glue together
solutions to the Seiberg--Witten equation for $-K_M$ on $M$
and reducible solutions on $N$ for the unique
$Spin^c$--structure with auxiliary line bundle
$e_1-\sum_{i\neq 1} e_i$, as in the proof
of Proposition 2 in~\cite{KMT}. This gives a Seiberg--Witten
invariant of $X$ which is equal (up to sign) to the
Seiberg--Witten invariant of $M$ for $-K_M$, which is
non--zero.

This implies that $L=-K_M+e_1-\sum_{i\neq 1} e_i$
has selfintersection number $=K_X^2$
because for symplectic manifolds all the non--trivial
Seiberg--Witten invariants come from zero--dimensional
moduli spaces, see~\cite{swg}. Thus, $a_i = \pm 1$
for all $i\in\{ 1,\ldots ,n\}$. Without loss of generality
we may assume $a_i =1$ for all $i$.

The line bundle $L$ is obtained from $-K_X$ by twisting
with $e_1$. Thus, by Taubes's main theorem in~\cite{swg},
the non--triviality of the Seiberg--Witten invariant of
$X$ with respect to $L$ implies that $e_1$ can be
represented by a symplectically embedded $2$--sphere
in $X$. This contradicts the minimality of $X$.

We conclude that $N$ must be an integral homology
sphere.
This completes the proof of Theorem~\ref{t:sympl}.

\begin{rem}
Gompf~\cite{gompf} has shown that all finitely presentable
groups occur as fundamental groups of minimal symplectic
$4$--manifolds, and conjecturally all these manifolds
are irreducible. As was the case in~\cite{IMRN,KMT},
our arguments do not give an optimal result because we
cannot deal with fundamental groups without non--trivial
finite quotients. With regard to Theorem~\ref{t:sympl},
note that there are such groups which occur as fundamental
groups of integral homology $4$--spheres. Let $G$ be
the Higman $4$--group, an infinite group without
non--trivial finite quotients, which has a presentation
with $4$ generators and $4$ relations. Doing surgery
on $4(S^1\times S^3)$ according to the relations
produces an integral homology sphere with fundamental
group $G$.
\end{rem}

\begin{rem}
In another direction, the assumption $b_2^+ (X)>1$ can
probably be removed from Theorem~\ref{t:sympl} and
Corollary~\ref{c:minirr}. To do this one needs to understand
how the neck--stretching in the proof of Theorem~\ref{t:sympl}
and the perturbations in Taubes's arguments~\cite{taubes,swg}
interact with the chamber structure of the Seiberg--Witten
invariants for manifolds with $b_2^+ =1$. We will return to
this question in the future.

However, some results about the case when $b_2^+ =1$
can be deduced from Theorem~\ref{t:sympl}. For example,
all manifolds with non--tivial finite fundamental groups
are dealt with by the following:
\end{rem}

\begin{cor}\label{c:b=1}
Let $X$ be a minimal symplectic $4$--manifold with
$b_2^+ (X)=1$ and $b_1 (X)\leq 1$. If $\pi_1 (X)$ is a
non--trivial residually finite group, then $X$ is irreducible.
\end{cor}
\begin{pf}
Suppose $X\cong M\# N$. We may assume that $N$ has negative
definite intersection form and its fundamental group has
no non--trivial finite quotient. Residual finiteness then
implies that $N$ is simply connected, and
$\pi_1 (M)\cong\pi_1 (X)$. By assumption, $X$ has a
finite cover $\overline{X}$ of degree $d>1$ which is diffeomorphic
to $\overline{M}\# dN$, where $\overline{M}$ is a $d$--fold
cover of $M$. But $\overline{X}$ is minimal symplectic
because $X$ is, and so
Corollary~\ref{c:minirr} implies that $N$ is a homotopy
sphere whenever $b_2^+ (\overline{X})>1$.

If $b_1 (X)=0$, the multiplicativity of the Euler
characteristic and of the signature imply $b_2^+
(\overline{X})\geq 3$. If $b_1 (X)=1$, we can take
$d\geq 3$ to obtain $b_2^+ (\overline{X})\geq 2$.
\end{pf}

\medskip\noindent
{\sl Acknowledgements:} I would like to thank Shuguang
Wang for correspondence and for his questions which
this paper answers, and Cliff Taubes for conversations
and correspondence.


\bibliographystyle{amsplain}

\bigskip

\end{document}